%
%%%%%%%%%%%%%%%%%%%%%%%%%%%%%%%%
%
% These are the macroes
%
%%%%%%%%%%%%%%%%%%%%%%%%%%%%%%%%
\newif\ifsect\newif\iffinal
\secttrue\finalfalse
\def\thm #1: #2{\medbreak\noindent{\bf #1:}\if(#2\thmp\else\thmn#2\fi}
\def\thmp #1) { (#1)\thmn{}}
\def\thmn#1#2\par{\enspace{\sl #1#2}\par
        \ifdim\lastskip<\medskipamount \removelastskip\penalty 55\medskip\fi}
\def\square{{\msam\char"03}}
\def\qedn{\thinspace\null\nobreak\hfill\square\par\medbreak}
\def\pf{\ifdim\lastskip<\smallskipamount \removelastskip\smallskip\fi
        \noindent{\sl Proof\/}:\enspace}
\def\itm#1{\item{\rm #1}\ignorespaces}

\def\bar#1{\overline{#1}}
%
%\input pdfcolor
%
%\input boxedeps
%\input boxedeps.cfg
%\HideDisplacementBoxes
%%
%\def\Figuraeps #1 (#2){\message{Figura #1}
%	\midinsert      
%	\centerline{\BoxedEPSF{#2.eps}}
%	\bigskip
%	\centerline{\bf Figure~#1}
%	\endinsert}
%%
%\def\Figurascaledeps #1 (#2)(#3){\message{Figura #1}
%	\midinsert      
%	\centerline{\BoxedEPSF{#2.eps scaled #3}}
%	\bigskip
%	\centerline{\bf Figure~#1}
%	\endinsert}
%
%
\newcount\parano
\newcount\eqnumbo
\newcount\thmno
\newcount\versiono
%\newcount\remno
\newbox\notaautore
\def\neweqt#1$${\xdef #1{(\number\parano.\number\eqnumbo)}
    \eqno #1$$
    \global \advance \eqnumbo by 1}
\def\newrem #1\par{\global \advance \thmno by 1
    \medbreak
{\bf Remark \the\parano.\the\thmno:}\enspace #1\par
\ifdim\lastskip<\medskipamount \removelastskip\penalty 55\medskip\fi}
\def\newex #1\par{\global \advance \thmno by 1
    \medbreak
{\it Example \the\parano.\the\thmno:}\enspace #1\par
\ifdim\lastskip<\medskipamount \removelastskip\penalty 55\medskip\fi}
\def\newthmt#1 #2: #3{ \global \advance \thmno by 1\xdef #2{\number\parano.\number\thmno}
    \medbreak\noindent
    {\bf #1 #2:}\if(#3\thmp\else\thmn#3\fi}
\def\neweqf#1$${\xdef #1{(\number\eqnumbo)}
    \eqno #1$$
    \global \advance \eqnumbo by 1}
\def\newthmf#1 #2: #3{    \global \advance \thmno by 1\xdef #2{\number\thmno}
    \medbreak\noindent
    {\bf #1 #2:}\if(#3\thmp\else\thmn#3\fi}
\def\forclose#1{\hfil\llap{$#1$}\hfilneg}
\def\newforclose#1{
	\ifsect\xdef #1{(\number\parano.\number\eqnumbo)}\else
	\xdef #1{(\number\eqnumbo)}\fi
	\hfil\llap{$#1$}\hfilneg
	\global \advance \eqnumbo by 1
	\iffinal\else\rsimb#1\fi}
\def\forevery#1#2$${\displaylines{\let\eqno=\forclose
        \let\neweq=\newforclose\hfilneg\rlap{$\qquad\quad\forall#1$}\hfil#2\cr}$$}
\def\noNota #1\par{}
\def\today{\ifcase\month\or
   January\or February\or March\or April\or May\or June\or July\or August\or
   September\or October\or November\or December\fi
   \space\number\year}
\def\inizia{\ifsect\let\neweq=\neweqt\else\let\neweq=\neweqf\fi
\ifsect\let\newthm=\newthmt\else\let\newthm=\newthmf\fi}
\def\bititolo{\empty}
\gdef\begin #1 #2\par{\xdef\titolo{#2}
\ifsect\let\neweq=\neweqt\else\let\neweq=\neweqf\fi
\ifsect\let\newthm=\newthmt\else\let\newthm=\newthmf\fi
\iffinal\let\Nota=\noNota\fi
\centerline{\titlefont\titolo}
\if\bititolo\empty\else\medskip\centerline{\titlefont\bititolo}
\xdef\titolo{\titolo\ \bititolo}\fi
\bigskip
\centerline{\bigfont
\autore \ifvoid\notaautore\else\footnote{${}^1$}{\unhbox\notaautore}\fi}
\bigskip\if\istituto!\centerline{\today}\else
\centerline{\istituto}
%\centerline{\indirizzo}
\centerline{\email}
\medskip
\centerline{#1\ \anno}\fi
\bigskip\bigskip
\ifsect\else\global\thmno=1\global\eqnumbo=1\fi}
\def\anno{2014}
\def\raggedleft{\leftskip2cm plus1fill \spaceskip.3333em \xspaceskip.5em
\parindent=0pt\relax}
\font\titlefont=cmssbx10 scaled \magstep1
\font\bigfont=cmr12
\font\eightrm=cmr8

\font\bbr=msbm10
\font\sbbr=msbm7
\font\ssbbr=msbm5
\font\msam=msam10

\font\bfm=cmmib10

\nopagenumbers
\binoppenalty=10000
\relpenalty=10000
\newfam\amsfam
\textfont\amsfam=\bbr \scriptfont\amsfam=\sbbr \scriptscriptfont\amsfam=\ssbbr
\newfam\boldifam
\textfont\boldifam=\bfm
\let\de=\partial

\let\phe=\varphi
\def\Hol{\mathop{\rm Hol}\nolimits}

\def\Re{\mathop{\rm Re}\nolimits}

\def\id{\mathop{\rm id}\nolimits}

\def\bigoperp{\mathop{\hbox{$\bigcirc\kern-11.8pt\perp$}}\limits}
\def\Klim{\mathop{\hbox{$K$-$\lim$}}\limits}
\def\rKlim{\mathop{\hbox{$K'$-$\lim$}}\limits}

\mathchardef\void="083F
\mathchardef\ellb="0960
\mathchardef\taub="091C
\def\C{{\mathchoice{\hbox{\bbr C}}{\hbox{\bbr C}}{\hbox{\sbbr C}}
{\hbox{\sbbr C}}}}
\def\R{{\mathchoice{\hbox{\bbr R}}{\hbox{\bbr R}}{\hbox{\sbbr R}}
{\hbox{\sbbr R}}}}

\newcount\notitle
\notitle=1
\headline={\ifodd\pageno\rhead\else\lhead\fi}
\def\rhead{\ifnum\pageno=\notitle\iffinal\hfill\else\hfill\tt Version
\the\versiono; \the\day/\the\month/\the\year\fi\else\hfill\eightrm\titolo\hfill
\folio\fi}
\def\lhead{\ifnum\pageno=\notitle\hfill\else\eightrm\folio\hfill\autore\hfill
\fi}
\newbox\bibliobox
\def\setref #1{\setbox\bibliobox=\hbox{[#1]\enspace}
    \parindent=\wd\bibliobox}
\def\biblap#1{\noindent\hang\rlap{[#1]\enspace}\indent\ignorespaces}
\def\art#1 #2: #3! #4! #5 #6 #7-#8 \par{\biblap{#1}#2: {\sl #3\/}.
    #4 {\bf #5} (#6)\if.#7\else, \hbox{#7--#8}\fi.\par}
\def\book#1 #2: #3! #4 \par{\biblap{#1}#2: {\bf #3.} #4.\par}
\def\coll#1 #2: #3! #4! #5 \par{\biblap{#1}#2: {\sl #3\/}. In {\bf #4,}
#5.\par}
\def\pre#1 #2: #3! #4! #5 \par{\biblap{#1}#2: {\sl #3\/}. #4, #5.\par}
\def\raggedleft{\leftskip2cm plus1fill \spaceskip.3333em \xspaceskip.5em
\parindent=0pt\relax}
\def\Nota #1\par{\medbreak\begingroup\Bittersweet\raggedleft
#1\par\endgroup\Black
\ifdim\lastskip<\medskipamount \removelastskip\penalty 55\medskip\fi}
%

%
%\newcount\defno
\def\smallsect #1. #2\par{\bigbreak\centerline{{\bf #1.}\enspace{\bf #2}}\sm\par
    \global\parano=#1\global\eqnumbo=1\global\thmno=0%\global\defno=0\global\remno=0
    \nobreak\smallskip\nobreak\noindent\message{#2}}
\def\newdef #1\par{\global \advance \thmno by 1
    \medbreak
{\bf Definition \the\parano.\the\thmno:}\enspace #1\par
\ifdim\lastskip<\medskipamount \removelastskip\penalty 55\medskip\fi}
%\finalfalse
%\versiono=12
\finaltrue

\magnification\magstephalf

%%%% Abbreviazioni

\let\bi=\bigskip

\let\sm=\smallskip

\def\rKlim{\mathop{\hbox{\rm $K'$-lim}}\limits}

\def\autore{Jasmin Raissy\footnote{}{\eightrm 2010 Mathematics Subject Classification: Primary 37L05; Secondary 32A40, 32H50, 20M20.\hfill\break\indent Keywords: infinitesimal generators, semigroups of holomorphic mappings, Julia-Wolff-Carath\'eodory theorem, boundary behaviour.}\footnote{}{\eightrm Partially supported by the FIRB2012 grant ``Differential Geometry and Geometric Function Theory'', and by the ANR project LAMBDA, ANR-13-BS01-0002.}}
\def\istituto{\vbox{\hfill Institut de Math\'ematiques de Toulouse; UMR5219,
Universit\'e de Toulouse; CNRS,\hfill\break\null\hfill
UPS IMT, F-31062 Toulouse Cedex 9, France.
E-mail: jraissy@math.univ-toulouse.fr\hfill\null}}
\def\email{}

\begin {November} The Julia-Wolff-Carath\'eodory theorem and its generalizations

This note is a short introduction to the Julia-Wolff-Carath\'eodory theorem, and its generalizations in several complex variables, up to very recent results for infinitesimal generators of semigroups.

%. A very precise and systematic presentation, providing clear proofs, of various aspects of the problem of generalization of the classical Julia-Wolff-Carath\'eodory theorem to several complex variables, and updated until 2004, can be found in [A6]. 

\smallsect 1.  The classical Julia-Wolff-Carath\'eodory theorem

One of the classical result in one-dimensional complex analysis is Fatou's theorem:

\newthm Theorem \zFatou: (Fatou [Fa]) Let $f\colon\Delta\to\Delta$ be a holomorphic self-map of the unit disk $\Delta\subset \C$. Then~$f$ admits non-tangential limit at almost every point of $\partial\Delta$.

This result however does not give precise information about the behavior
at a specific point~$\sigma$ of the boundary.
Of course, to obtain a more precise statement in this case some hypotheses
on~$f$ are needed. In fact, as it was found by Julia ([Ju1]) in  
1920, the right hypothesis is to assume that $f(\zeta)$
approaches the boundary of~$\Delta$ at least as fast as~$\zeta$, in a weak sense.
More precisely, we have the classical {\it Julia's lemma:}

\newthm Theorem \zJulia: (Julia [Ju1])
Let $f\colon\Delta\to\Delta$ be a
bounded holomorphic function such that
$$
\liminf_{\zeta\to\sigma}{1-|f(\zeta)|\over1-|\zeta|}=\alpha<+\infty
\neweq\eqzalp
$$
for some $\sigma\in\de\Delta$. Then $f$ has non-tangential limit $\tau\in\de
\Delta$ at~$\sigma$. Moreover, for all $\zeta\in\Delta$ one has
$$
{|\tau-f(\zeta)|^2\over1-|f(\zeta)|^2}\le\alpha\,{|\sigma-\zeta|^2\over1-
	|\zeta|^2}\;.
\neweq\eqzJ
$$

The latter statement admits an interesting geometrical interpretation. The {\it
horo\-cycle} $E(\sigma,R)$ contained in $\Delta$ of {\it center}~$\sigma\in\de\Delta$ and
{\it radius}~$R>0$ is the set
$$
E(\sigma,R)=\left\{\zeta\in\Delta\biggm|{|\sigma-\zeta|^2\over1-|\zeta|^2}
	<R\right\}\;.
$$
Geometrically, $E(\sigma,R)$ is an euclidean disk of radius $R/(R+1)$
internally tangent to~$\de\Delta$ at~$\sigma$. Therefore~\eqzJ\ becomes
$f\bigl(E(\sigma,R)\bigr)\subseteq E(\tau,\alpha R)$ for all~$R>0$, and the
existence of the non-tangential limit more or less follows from~\eqzJ\ and
from the fact that horocycles touch the boundary in exactly one point.

A horocycle can be thought of as the limit of Poincar\'e disks of fixed
euclidean radius and centers going to the boundary; so it makes sense to think of horocycles as Poincar\'e disks centered at the boundary, and of Julia's
lemma as a Schwarz-Pick lemma at the boundary. This suggests that $\alpha$
might be considered as a sort of dilation coefficient: $f$ expands horocycles
by a ratio of~$\alpha$. If $\sigma$ were an internal point and $E(\sigma,R)$ an
infinitesimal euclidean disk actually centered at~$\sigma$, one then would be tempted
to say that~$\alpha$ is (the absolute value of) the derivative of~$f$ at~$\sigma$.
This is exactly the content of the classical {\it  Julia-Wolff-Carath\'eodory
theorem:}

\newthm Theorem \zJWC: (Julia-Wolff-Carath\'eodory) Let $f\colon\Delta\to\Delta$ be a
bounded holomorphic function such that
$$
\liminf_{\zeta\to\sigma}{1-|f(\zeta)|\over1-|\zeta|}=\alpha<+\infty
$$
for some $\sigma\in\de\Delta$, and let $\tau\in\de\Delta$ be the non-tangential
limit of~$f$ at~$\sigma$. Then both the incremental ratio $\bigl(\tau-f(\zeta)
\bigr)\big/(\sigma-\zeta)$ and the derivative~$f'(\zeta)$ have 
non-tangential limit~$\alpha\bar\sigma\tau$ at~$\sigma$.

So condition \eqzalp\ forces the existence of the non-tangential limit of
both~$f$ and its derivative at~$\sigma$. This is the result of the work of
several people: Julia~[Ju2], Wolff~[Wo], Carath\'eodory~[C], Landau and
Valiron~[L-V], R.~Nevanlinna~[N] and others. We refer, for example, to~[B] and~[A1] for proofs,
history and applications.

\bi

\smallsect 2. Generalizations to several variables

It was first remarked by Kor\'anyi and Stein ([Ko], [K-S], [St]) in extending Fatou's theorem to several complex variables, that the notion of non-tangential limit is not the right one to consider for domains in $\C^n$. In fact, it turns out that two notions are needed, and to introduce them it is useful to investigate the notion of non-tangential limit in the unit disk~$\Delta$. 

The non-tangential limit can be defined in two equivalent ways. A
function $f\colon\Delta\to\C$ is said to have {\it non-tangential limit}~$L\in\C$
at~$\sigma\in\de\Delta$ if~$f\bigl(\gamma(t)\bigr)\to L$ as~$t\to 1^-$ for
every curve $\gamma\colon[0,1)\to\Delta$ such that $\gamma(t)$ converges to $\sigma$
non-tangentially as~$t\to1^-$. In $\C$, this is equivalent to having that~$f(\zeta)\to L$ as~$\zeta\to\sigma$ staying inside any {\it
Stolz region}~$K(\sigma, M)$ of {\it vertex}~$\sigma$ and~{\it
amplitude}~$M>1$, where
$$
K(\sigma,M)=\left\{\zeta\in\Delta\biggm|{|\sigma-\zeta|\over1-|\zeta|}<M
	\right\}\;,
$$
since Stolz regions are angle-shaped nearby the vertex~$\sigma$, and the angle
is going to~$\pi$ as~$M\to+\infty$. These two approaches lead to different notions in several variables.

In the unit ball $B^n\subset\C^n$ the natural generalization of a Stolz region is the {\it Kor\'anyi region} $K(p, M)$ of {\it vertex}~$p\in\de B^n$ and~{\it
amplitude}~$M>1$ given by
$$
K(p,M)=\left\{z\in B^n\biggm|{|1-\langle z,p\rangle|\over1-\|z\|}<M
	\right\}\;,
$$
where $\|\cdot\|$ denote the euclidean norm and $\langle\cdot\,,\cdot\rangle$ the canonical hermitian product. Then a function $f\colon B^n\to\C$ has {\it $K$-limit} (or {\it admissible limit})~$L\in\C$ at~$p\in\de B^n$, and we write
$$
\Klim_{z\to p} f(z)
$$
if $f(z)\to L$ as~$z\to p$ staying inside any Kor\'anyi region~$K(\sigma, M)$. A Kor\'anyi region $K(p,M)$ approaches the boundary non-tangentially along the normal direction at~$p$ but tangentially along the complex tangential directions at~$p$. Therefore, having $K$-limit is stronger than having non-tangential limit. However, as first noticed by Kor\'anyi and Stein, for holomorphic functions of several complex variables one is often able to prove the existence of $K$-limits.  For instance, the best generalization of Julia's lemma to $B^n$ is the following result (proved by Herv\'e [H] in terms of non-tangential limits and by Rudin [R] in general):

\newthm Theorem \zJdue: (Rudin [R]) Let $f\colon B^n\to B^m$ be a holomorphic map such that 
$$
\liminf_{z\to p}{1-\|f(z)\|\over1-\|z\|}=\alpha<+\infty\;,
$$
for some $p\in\de B^n$. Then $f$ admits $K$-limit $q\in\de B^m$ at~$p$, and furthermore for all $z\in B^n$ one has
$$
{|1-\langle f(z),q\rangle|^2\over 1-\|f(z)\|^2}\le \alpha {|1-\langle z,p\rangle|^2\over 1-\|z\|^2}\;.
$$

To define Kor\'anyi regions for more general domains in $\C^n$ than the unit ball, we need to briefly recall the definition of the Kobayashi distance (we refer, e.g., to [A1], [JP] and [Ko] for details and much more).  We denote by~$k_\Delta$ the Poincar\'e distance on the unit disk~$\Delta\subset\C$. Given $X$ a complex manifold,
the {\sl Lempert function}~$\delta_X\colon X\times X\to\R^+$ of~$X$ is defined as
$$
\delta_X(z,w)=\inf\{k_\Delta(\zeta,\eta)\mid
\hbox{$\exists\phi\colon\Delta\to X$ holomorphic, with $\phi(\zeta)=z$ and
$\phi(\eta)=w$}\}
$$
for all $z$, $w\in X$. The {\sl Kobayashi pseudodistance}~$k_X\colon X\times X\to\R^+$ of~$X$
is then defined as the largest pseudodistance on~$X$ bounded above by~$\delta_X$. The manifold $X$ is called {\sl (Kobayashi) hyperbolic} if $k_X$ is indeed a distance; $X$ is called {\sl complete hyperbolic} if $k_X$ is a complete distance. 

The main property of the Kobayashi (pseudo)distance is that it is contracted by holomorphic maps: if $f\colon X\to Y$ is a holomorphic map then
$$
\forevery{z,w\in X} k_Y\bigl(f(z),f(w)\bigr)\le k_X(z,w)\;.
$$
In particular, the Kobayashi distance is invariant under biholomorphisms.

It is easy to see that the Kobayashi distance of the unit disk coincides with the Poincar\'e distance. Furthermore, the Kobayashi distance of the unit ball $B^n\subset\C^n$ coincides with the
Bergman distance (see, e.g., [A1, Corollary~2.3.6]); and the Kobayashi distance of a bounded
convex domain coincides with the Lempert function (see, e.g., [A1, Proposition~2.3.44]). Moreover, the Kobayashi distance of a bounded convex domain $D$ is complete ([A1, Proposition~2.3.45]), and thus for each $p\in D$ we have that $k_D(p,z)\to+\infty$ if and only if $z$ tends to the boundary $\de D$. 

Using the Kobayashi intrinsic distance we obtain the natural generalization to complete hyperbolic domains of Kor\'anyi regions of the balls.

Let $D\subset\subset \C^n$ be a complete hyperbolic domain and denote by $k_D$ its Kobayashi distance.
A {\sl $K$-region} of {\sl vertex $x\in\partial D$, amplitude $M>1$}, and {\sl pole $z_0\in D$} is the set
$$
K_{D,z_0}(x,M) = \left\{z\in D \mid \limsup_{w\to x}\left[k_D(z,w) - k_D(z_0,w)\right] + k_D(z_0, z)< \log M\right\}.
$$
This definition clearly depends on the pole $z_0$. However, this dependence is not too relevant since changing the pole corresponds to shifting amplitudes. Moreover, it is elementary to check that in the unit ball $K$-regions coincide with Kor\'anyi regions, $K_{B^n,0}(x,M) = K(x,M)$. Therefore $K$-regions are a natural generalization of Kor\'anyi regions allowing us to generalize the notion of $K$-limit. A function $f\colon D\to\C^m$ has {\it $K$-limit}~$L$ at~$x\in\de D$ if $f(z)\to L$ as~$z\to p$ staying inside any $K$-region of vertex $x$. The best generalization of Julia's lemma in this setting is then the following, due to Abate:

\newthm Theorem \zJdueAba: (Abate [A2]) Let $D\subset\subset \C^n$ be a complete hyperbolic domain and let $z_0\in D$. Let $f\colon D\to \Delta$ be a holomorphic function and let $x\in\partial D$ be such that 
$$
\liminf_{z\to x}\left[k_D(z_0,z) - k_\Delta(0, f(z))\right] <+\infty\;.
$$
Then $f$ admits $K$-limit $\tau\in\de D$ at~$x$.

\sm

In order to obtain a complete generalization of the 
Julia-Wolff-Carath\'eodory for~$B^n$, Rudin discovered that he needed a different notion of limit, still stronger than non-tangential limit but weaker than $K$-limit. This notion is closely related to the other characterization of non-tangential limit in one variable we mentioned at the beginning of this section.

A crucial one-variable result relating limits along curves and non-tangential limits is {\it Lindel\"of's theorem.} Given $\sigma\in\de\Delta$, a {\it $\sigma$-curve} is a continuous curve $\gamma\colon[0,1)\to\Delta$ such that $\gamma(t)\to\sigma$ as $t\to 1^-$. Then Lindel\"of [Li] proved that if a bounded holomorphic function $f\colon\Delta\to\C$ admits limit~$L\in\C$ along a given $\sigma$-curve then it admits limit $L$ along all non-tangential $\sigma$-curves --- and thus it has non-tangential limit~$L$ at~$\sigma$. 

In generalizing this result to several complex variables, \v Cirka [\v C] realized that for a bounded holomorphic function the existence of the limit along a (suitable) $p$-curve (where $p\in\de B^n$) implies not only the existence of the non-tangential limit, but also the existence of the limit along any curve belonging to a larger class of curves, including some tangential ones --- but it does not in general imply the existence of the $K$-limit. To describe the version (due to Rudin [R]) of \v Cirka's result we shall state in this survey, let us introduce a bit of terminology. 

Let $p\in\de B^n$. As before, a {\it $p$-curve} is a continuous curve $\gamma\colon[0,1)\to B^n$ such that
$\gamma(t)\to p$ as $t\to 1^-$. A $p$-curve is {\it special} if
$$
\lim_{t\to 1^-}{\|\gamma(t)-\langle\gamma(t),p\rangle p\|^2\over 1-|\langle\gamma(t),p\rangle|^2}=0\;;
\neweq\eqzspec
$$
and, given $M>1$, it is {\it $M$-restricted} if
$$
{|1-\langle\gamma(t),p\rangle|\over 1-|\langle\gamma(t),p\rangle|}<M
$$
for all $t\in[0,1)$. We also say that $\gamma$ is {\it restricted} if it is $M$-restricted for some~$M>1$. In other words, $\gamma$ is restricted if and only if $t\mapsto\langle\gamma(t),p\rangle$ goes to~1 non-tangentially in~$\Delta$. 

It is not difficult to see that non-tangential curves are special and restricted; on the other hand,
a special restricted curve approaches the boundary non-tangentially along the normal direction, but it can approach the boundary tangentially along complex tangential directions. Furthermore,
a special $M$-restricted $p$-curve is eventually contained in any $K(p,M')$ with $M'>M$,
and conversely a special $p$-curve eventually contained in $K(p,M)$ is $M$-restricted. 
However, $K(p,M)$ can contain $p$-curves that are restricted but not special: for these curves the limit in \eqzspec\ might be a strictly positive number. 

With these definitions in place, we shall say that a function $f\colon B^n\to\C$ has
{\it restricted $K$-limit} (or {\it hypoadmissible limit})~$L\in\C$ at~$p\in\de B^n$, and we shall write
$$
\rKlim_{z\to p} f(z)=L\;,
$$
if $f\bigl(\gamma(t)\bigr)\to L$ as~$t\to 1^-$ for any special restricted $p$-curve $\gamma\colon[0,1)\to B^n$. It is clear that the existence of the $K$-limit implies the existence of the restricted $K$-limit, that in turns implies the existence of the non-tangential limit; but none of these implications can be reversed (see, e.g., [R] for examples in the ball). 

Finally, we say that a function $f\colon B^n\to\C$ is
{\it $K$-bounded} at~$p\in\de B^n$ if it is bounded in any Kor\'anyi region $K(p,M)$, where the bound can depend on~$M>1$. Then Rudin's version of \v Cirka's generalization of Lindel\"of's theorem is the following:

\newthm Theorem \zCR: (Rudin [R])
Let $f\colon B^n\to\C$ be a holomorphic function $K$-bounded at $p\in\de B^n$. Assume there is a special restricted $p$-curve $\gamma^o\colon[0,1)\to B^n$ such that $f\bigl(\gamma^o(t)\bigr)\to L\in\C$ as~$t\to 1^-$. Then $f$ has restricted $K$-limit $L$ at~$p$.

As before, it is possible to generalize this approach to a domain $D\subset\C^n$ different from the ball. A very precise and systematic presentation, providing clear proofs, details and examples, of various aspects of the problem of generalization of the classical Julia-Wolff-Carath\'eodory theorem to domains in several complex variables, and updated until 2004, can be found in [A6]. 

For the sake of simplicity we state here only the definitions needed to state Abate's version of Lindel\"of's theorem in this setting. 
%The interested reader can find details and examples in [A6]. 
Given a point $x\in\de D$, a {\it $x$-curve} is again a continuous curve $\gamma\colon [0,1)\to D$ so that $\lim_{t\to 1^-} \gamma(t) = x$. A {\it projection device} at $x\in\de D$ is the data of: a neighbourhood $U$ of $x$ in $\C^n$, a holomorphic embedded disk $\phe_x\colon \Delta\to D\cap U$, such that $\lim_{\zeta\to 1} \phe_x(\zeta) = x$, a family $\cal P$ of $x$-curves in $D\cap U$, and a device associating to every $x$-curve $\gamma\in\cal P$ a 1-curve $\tilde \gamma_x$ in $\Delta$, or equivalently a $x$-curve $\gamma_x = \phe_x\circ \tilde \gamma_x$ in $\phe_x(\Delta)$. 
If $D$ is equipped with a projection device at $x\in\de D$, then a curve $\gamma\in \cal P$ is {\it special} if $\lim_{t\to 1^-} k_{D\cap U}(\gamma(t), \gamma_x(t)) = 0$,
and it is {\it restricted} if $\gamma_x$ is a non-tangential $1$-curve in $\Delta$. A function $f\colon D\to C$ has {\it restricted $K$-limit $L\in\C$ at $x$} if $\lim_{t\to 1^-} f(\gamma(t)) = L$ for all special restricted $x$-curves.  
A projection device is {\it good} if: for any $M>1$ there is a $M'>1$ so that $\phe_x(K(1,M)) \subset K_{D\cap U, z_0}(x, M')$, and for any special restricted $x$-curve $\gamma$ there exists $M_1 = M_1(\gamma)$ such that $\lim_{t\to 1^-} k_{K_{D\cap U, z_0}(x, M_1)}(\gamma(t), \gamma_x(t)) = 0$. Good projection devices exist, and several examples can be found for example in [A6].
Finally, we say that a function $f\colon D\to\C$ is
{\it $K$-bounded} at~$p\in\de B^n$ if it is bounded in any $K$-region $K_{D, z_0}(x,M)$, where the bound can depend on~$M>1$. 

With these definitions we can state the generalization of Lindel\"of principle given by Abate.

\newthm Theorem \zCRAba: (Abate [A6])
Let $D\subset\C^n$ be a domain equipped with a good projection device at $x\in \de D$. Let $f\colon D\to\Delta$ be a holomorphic function $K$-bounded at $x$. Assume there is a special restricted $x$-curve $\gamma^o\colon[0,1)\to D$ such that $f\bigl(\gamma^o(t)\bigr)\to L\in\C$ as~$t\to 1^-$. Then $f$ has restricted $K$-limit $L$ at~$x$.

We can now deal with the generalization of the Julia-Wolff-Carath\'eodory theorem to several complex variables. With respect to the one-dimensional case there is an obvious difference:
instead of only one derivative one has to deal with a whole (Jacobian) matrix of
them, and there is no reason they should all behave in the same way. And
indeed they do not, as shown in
Rudin's version of the Julia-Wolff-Carath\'eodory theorem for the unit ball:

\newthm Theorem \zJWCR: (Rudin [R])
Let $f\colon B^n\to B^m$ be a holomorphic map such that 
$$
\liminf_{z\to p}{1-\|f(z)\|\over1-\|z\|}=\alpha<+\infty\;,
$$
for some $p\in\de B^n$. Then $f$ admits 
$K$-limit $q\in\de B^m$ at $p$. Furthermore, if we set
$f_q(z)=\bigl\langle f(z),p
\bigr\rangle q$ and
denote by~$df_z$ the differential of~$f$ at~$z$, we have:
\item{\rm (i)} the function
$\bigl[1-\bigl\langle f(z),q\bigr\rangle\bigr]\big/[1-\langle z,p\rangle]$ is $K$-bounded and has restricted $K$-limit~$\alpha$ at~$p$;
\item{\rm (ii)} the map $[f(z)-f_q(z)]/[1-\langle z,p\rangle]^{1/2}$ is $K$-bounded and has restricted
$K$-limit~$O$ at~$p$;
\item{\rm (iii)} the function $\bigl\langle df_z(p),q\bigr\rangle$ is $K$-bounded and has restricted
$K$-limit~$\alpha$ at~$p$;
\item{\rm (iv)} the map $[1-\langle z,p\rangle]^{1/2}d(f-f_q)_z(p)$ is $K$-bounded and has restricted
$K$-limit~$O$ at~$p$;
\item{\rm (v)} if $v$ is any vector orthogonal to~$p$, the function
$\bigl\langle df_z (v),q\bigr\rangle\big/[1-\langle z,p\rangle]^{1/2}$ 
is $K$-bounded and has restricted $K$-limit~$0$ at~$p$;
\item{\rm (vi)} if $v$ is any vector orthogonal to~$p$, the map
$d(f-f_q)_z (v)$ is $K$-bounded at~$p$.

In the last twenty years this theorem (as well as Theorems~\zJdue\ and~\zCR) has been extended to domains much more general than the unit ball: for instance, strongly pseudoconvex domains [A1, 2, 3], convex domains of finite type [AT], and polydisks [A5] and [AMY], (see also [A6] and references therein). 

We end this section with the general version of the Julia-Wolff-Carath\'odory theorem obtained by Abate in [A6] for a complete hyperbolic domain $D$ in $\C^n$. To formulate it, we need to introduce a couple more definitions. A projection device at $x\in\de D$ is {\it geometrical} if there is a holomorphic function $\tilde p_x\colon D\cap U\to \Delta$ such that $\tilde p_x\circ \phe_x = \id_\Delta$ and $\tilde\gamma_x = \tilde p_x\circ \gamma$ for all $\gamma \in \cal P$. A geometrical projection device at $x$ is {\it bounded} if $d(z,\de D)/|1-\tilde p_x(z)|$ is bounded in $D\cap U$, and $|1-\tilde p_x(z)|/d(z,\de D)$ is $K$-bounded in $D\cap U$. The statement is then the following, where $\kappa_D$ denotes the Kobayashi metric.
  
\newthm Theorem \zJWCAba: (Abate [A6])
Let $D\subset\C^n$ be a complete hyperbolic domain equipped with a bounded geometrical projection device at $x\in\de D$. Let $f\colon D\to \Delta$ be a holomorphic function such that 
$$
\liminf_{z\to x}\left[k_D(z_0,z) - k_\Delta(0, f(z))\right] = {1\over 2}\log\beta<+\infty\;.
$$
Then for every $v\in\C^n$ and every $s\ge0$ such that $d(z, \de D)^s \kappa_D(z;v)$ is $K$-bounded at $x$ the function
$$
d(z, \de D)^{s-1} {\de f\over\de v}
\neweq\stella
$$
is $K$-bounded at $x$. Moreover, if $s>\inf\{s\ge0\mid d(z, \de D)^s \kappa_D(z;v)~\hbox{is}~K\hbox{-bounded at}~x\}$, then \stella\ has vanishing $K$-limit at $x$.

Depending on more specific properties of the projection device, it is indeed possible to deduce the existence of restricted $K$-limits, see [A6, Section 5].

Further generalizations of Julia-Wolff-Carath\'eodory theorem have been obtained in infinite-dimen\-sio\-nal Banach and Hilbert spaces, and we refer to [EHRS], [ELRS], [ERS], [F], [MM], [SW], [W\l 1, 2, 3], [Z], and references therein.

\smallsect 3. Julia-Wolff-Carath\'eodory theorem for infinitesimal generators

We conclude this survey focusing on a different kind of generalization in several complex variables:  infinitesimal generators of one-parameter semigroups of holomorphic self-maps of~$B^n$.

We consider $\Hol(B^n,B^n)$, the space of holomorphic self-maps of~$B^n$, endowed with the usual compact-open topology. A {\it one-parameter semigroup} of holomorphic self-maps of~$B^n$ is a continuous semigroup homomorphism $\Phi\colon\R^+\to\Hol(B^n,B^n)$. In other words, writing $\phe_t$ instead of~$\Phi(t)$, we have $\phe_0=\id_{B^n}$, the map $t\mapsto\phe_t$ is continuous, and the semigroup property
$\phe_t\circ\phe_s=\phe_{t+s}$
holds. An introduction to the theory of one-parameter semigroups of holomorphic maps can be found in [A1], [RS2] or [S].

One-parameter semigroups can be seen as the flow of a vector field (see, e.g., [A4]). Given a 
one-parameter semigroup $\Phi$, it is possible to prove that there exists a holomorphic map 
$G\colon B^n\to\C^n$, the {\it infinitesimal generator} of the semigroup, such that
$$
{\de\Phi\over\de t}=G\circ\Phi\;.
\neweq\eqzinfgen
$$
It should be kept in mind, when reading the literature on this subject, that in some papers (e.g., in [ERS] and [RS1]) there is a change of sign with respect to our definition, due to the fact that the infinitesimal generator is defined there as the solution of the equation
$$
{\de\Phi\over\de t}+G\circ\Phi=O\;.
$$

A Julia's lemma for infinitesimal generators was proved by Elin, Reich and Shoikhet in [ERS] in 2008, assuming that the radial limit of the generator at a point $p\in\de B^n$ vanishes:

\newthm Theorem \ERS: ([ERS, Theorem p. 403]) Let $G\colon B^n\to\C^n$ be the infinitesimal generator on $B^n$ of a one-parameter semigroup $\Phi=\{\phe_t\}$, and let $p\in\de B^n$ be such that
$$
\lim_{t\to 1^-}G(tp)=O\;.
\neweq\eqzbrnp
$$
Then the following assertions are equivalent:
\smallskip
\itm{(I)} $\alpha=\liminf_{t\to 1^-} \Re{\langle G(tp),p\rangle\over t-1}<+\infty$;
\itm{(II)} $\beta= 2\sup_{z\in B^n}\Re \left[{\langle G(z),z\rangle\over 1-\|z\|^2}-{\langle G(z),p\rangle\over 1-\langle z,p\rangle}
\right]<+\infty$;
\itm{(III)} there exists $\gamma\in\R$ such that for all $z\in B^n$ we have ${|1-\langle\phe_t(z),p\rangle|^2\over 1-\|\phe_t(z)\|^2}\le e^{\gamma t}{|1-\langle z,p\rangle|^2\over 1-\|z\|^2}$.
\smallskip
\noindent Furthermore, if any of these assertions holds then $\alpha=\beta =\inf\gamma$, and
we have
$$
\lim_{t\to 1^-}{\langle G(tp),p\rangle\over t-1}=\beta\;.
\neweq\eqzbeta
$$

If \eqzbrnp\ and any (whence all) of the equivalent conditions (I)--(III) holds,~$p\in\de B^n$ is called a {\it boundary regular null point} of~$G$ with {\it dilation}~$\beta\in\R$. 

This result suggested that a Julia-Wolff-Carath\'eodory theorem could hold 
for infinitesimal generators along the line of Rudin's Theorem~\zJWCR. A first partial generalization has been achieved by Bracci and Shoikhet in [BS]. In collaboration with Abate, in [AR] we have been able to give a full generalization of Julia-Wolff-Carath\'eodory theorem for infinitesimal generators, proving the following result.

\newthm Theorem \zJWCAR: ([AR]) Let $G\colon B^n\to\C^n$ be an infinitesimal generator on $B^n$ of a one-parameter semigroup, and let $p\in\de B^n$. Assume that
$$
{\langle G(z), p\rangle\over \langle z,p\rangle-1}\hbox{\qquad is $K$-bounded at $p$}
\neweq\eqzassdue
$$
and
$$
{G(z)-\langle G(z),p\rangle p\over(\langle z,p\rangle-1)^{\gamma}}\hbox{\qquad is $K$-bounded at $p$ for some $0<\gamma<1/2$.}
\neweq\eqzassdueb
$$
Then $p\in\de B^n$ is a boundary regular null point  for $G$. Furthermore, if $\beta$ is the dilation of $G$ at~$p$ then:
\itm{(i)} the function
$\langle G(z),p\rangle\big/(\langle z,p\rangle-1)$ (is $K$-bounded and) has restricted $K$-limit~$\beta$ at~$p$;
\itm{(ii)} if $v$ is a vector orthogonal to $p$, the function $\langle G(z),v\rangle/(\langle z,p\rangle-1)^{\gamma}$ is $K$-bounded and has restricted $K$-limit $0$ at~$p$;
\itm{(iii)} the function $\langle dG_z(p),p\rangle$ is $K$-bounded and has restricted $K$-limit $\beta$ at~$p$;
\itm{(iv)} if $v$ is a vector orthogonal to~$p$, the function $(\langle z,p\rangle-1)^{1-\gamma}\langle dG_z(p),v\rangle$ is $K$-bounded and has restricted $K$-limit~$0$ at~$p$;
\itm{(v)} if $v$ is a vector orthogonal to~$p$, the function $\langle dG_z(v),p\rangle\bigm/ (\langle z,p\rangle-1)^{\gamma}$ is $K$-bounded and has restricted $K$-limit $0$ at $p$.
\itm{(vi)} if $v_1$ and $v_2$ are vectors orthogonal to~$p$ the function $(\langle z,p\rangle-1)^{1/2-\gamma}\langle dG_z(v_1),v_2\rangle$ is $K$-bounded\break\indent at~$p$.

\sm\noindent{\sl Sketch of Proof of Theorem~\zJWCAR\/.\enspace} Statement (i) follows immediately from our hypotheses, thanks to Theorems~\zCR\ and~\ERS. Statement (iii) follows by standard arguments, and (iv) follows from (ii), again by standard arguments.

The main point is the proof of part (ii). 
By Theorem~\zCR, it suffices to compute the limit along a special restricted curve.
We use the curve
$$
\sigma(t)=t p + e^{-i\theta}\varepsilon(1-t)^{1-\gamma} v
$$
which is always restricted, and {\sl it is special if and only if $\gamma<1/2$}. We then plug (i) and this curve into Theorem~\ERS.(II), and we then let $\varepsilon\to 0^+$, using $\theta$ to get rid of the real part.

Statement (v) follows from (i), (ii) and by Theorem~\ERS\ using somewhat delicate arguments involving a curve of the form
$$
\gamma(t)=\bigl(t+ic(1-t)\bigr) p + \eta(t) v\;,
$$
where $1-t<|\eta(t)|^2<1-|t+ic(1-t)|^2$, and the argument of $\eta(t)$ is chosen suitably.
\qedn

A first difference with respect to Theorem~\zJWCR\ is that we have to assume \eqzassdue\ and \eqzassdueb\ as separate hypotheses, whereas they appear as part of Theorem~\zJWCR.(i) and (ii). Indeed, when dealing with holomorphic maps, conditions \eqzassdue\ and \eqzassdueb\ are a consequence of the equivalent of condition (I) in Theorem~\ERS, but in that setting the proof relies in the fact that there we have holomorphic {\sl self-maps} of the ball. 
In our context, \eqzassdueb\ is {\sl not} a consequence of Theorem~\ERS.(I), as shown in [AR, Example~1.2]; and it also seems that \eqzassdue\ is stronger than Theorem~\ERS.(I). 

A second difference is the exponent $\gamma <1/2$. Bracci and Shoikhet proved Theorem \zJWCAR\ with $\gamma=1/2$ but they couldn't prove the statements about restricted $K$-limits in cases (ii), (iv) and (v). This is due to an obstruction, which is not just a technical problem, but an inevitable feature of the theory. As mentioned in the sketch of the proof, in showing the existence of restricted $K$-limits, the curves one would like to use for obtaining the exponent $1/2$ in the statements are {\sl restricted but not special}, in the sense that the limit in~\eqzspec\ is a strictly positive (though finite) number. Actually the exponent $1/2$ might not be the right one to consider in the setting of infinitesimal generators, as shown in [AR, Example 1.2]. 

\sm An exact analogue of Theorem~\zJWCR\ with $\gamma=1/2$ can be recovered assuming a slightly stronger hypothesis on the infinitesimal generator.
Under assumptions \eqzassdue\ and \eqzassdueb\ we have
$$
{\bigl\langle G\bigl(\sigma(t)\bigr),p\bigr\rangle\over \langle \sigma(t),p\rangle-1} = \beta+o(1)
\neweq\eqJC
$$
as $t\to 1^-$ for any special restricted $p$-curve $\sigma\colon[0,1)\to B^n$.
Following ideas introduced in [ESY], [EKRS] and [EJ] in the context of the unit disk, $p$ is said to be a {\it H\"older boundary null point} if there is $\alpha>0$ such that 
$$
{\bigl\langle G\bigl(\sigma(t)\bigr),p\bigr\rangle\over \langle \sigma(t),p\rangle-1} = \beta+o\bigl((1-t)^\alpha\bigr)
\neweq\eqH
$$
for any special restricted $p$-curve $\sigma\colon[0,1)\to B^n$ such that $\langle\sigma(t),p\rangle\equiv t$. Using this notion we obtain the following result.

\newthm Theorem \zJWCARb: ([AR]) Let $G\colon B^n\to\C^n$ be the infinitesimal generator on $B^n$ of a one-parameter semigroup, and let $p\in\de B^n$. Assume that
$$
{\langle G(z), p\rangle\over \langle z,p\rangle-1}\quad\hbox{and}\quad
{G(z)-\langle G(z),p\rangle p\over(\langle z,p\rangle-1)^{1/2}} 
$$
are $K$-bounded at $p$, and that $p$ is a H\"older boundary null point. Then the statement of Theorem~\zJWCAR\ holds with $\gamma=1/2$.

Examples of infinitesimal generators with a H\"older boundary null point and satisfying the hypotheses of Theorem~\zJWCARb\ are provided in [AR].

\sm In a forthcoming paper in collaboration with Abate, we will deal with the generalization of Theorem \zJWCAR\ to strongly convex domains in $\C^n$.

\setref{EKRS}
\beginsection References

\book A1 M. Abate: Iteration theory of holomorphic maps on taut manifolds!
Me\-di\-ter\-ranean Press, Rende, 1989

\art A2 M. Abate: The Lindel\"of principle and the angular derivative in
strongly convex domains! J. Analyse Math.! 54 1990 189-228

\art A3 M. Abate: Angular derivatives in strongly pseudoconvex domains! Proc. Symp. 
Pure Math.! {52, \rm Part 2} 1991 23-40

\art A4 M. Abate: The infinitesimal generators of semigroups of holomorphic maps! Ann. Mat. Pura Appl.! 161 1992 167-180

\art A5 M. Abate: The Julia-Wolff-Carath\'eodory theorem in polydisks! J.
Analyse Math.! 74 1998 275-306 

\coll A6 M. Abate: Angular derivatives in several complex variables! Real methods in
complex and CR geometry! Eds. D. Zaitsev, G. Zampieri, Lect. Notes in Math.
1848, Springer, Berlin, 2004, pp. 1--47 

\art AR M. Abate, J. Raissy: A Julia-Wolff-Carath\'eodory theorem for infinitesimal generators in the unit ball! to appear in Trans. AMS,! {} 2014 1-17

\art AT M. Abate,  R. Tauraso: The Lindel\"of principle and angular derivatives
in convex domains of finite type! J. Austr. Math. Soc.! 73 
2002 221-250

\art AMY J. Agler, J.E. McCarthy, N.J. Young: A Carath\'eodory theorem for the bidisk via Hilbert space methods! Math. Ann.! 352 2012 581-624

\art BCD F. Bracci, M.D. Contreras, S. D\'\i az-Madrigal: Pluripotential theory, semigroups and boun\-dary behavior of infinitesimal generators in strongly convex domains! J. Eur. Math. Soc.! 12 2010 23-53

\art BS F. Bracci, D. Shoikhet: Boundary behavior of infinitesimal generators in the unit ball! Trans. Amer. Math. Soc.! 366 2014 1119-1140

\book B R.B. Burckel: An introduction to classical complex analysis! Academic Press, New York,1979

\art C C. Carath\'eodory: \"Uber die Winkelderivierten von beschr\"ankten
analytischen Funktionen! Sitzungsber. Preuss. Akad. Wiss. Berlin! {} 1929 39-54

\art {\v C} E.M. \v Cirka: The Lindel\"of and Fatou theorems in $\C^n$! Math.
USSR-Sb.! 21 1973 619-641

\art EHRS M. Elin, L.A. Harris, S. Reich, D. Shoikhet: Evolution equations and geometric function theory in $J^*$-algebras! J. Nonlinear Conv. Anal.! 3 2002 81-121

\art EKRS M. Elin, D. Khavinson, S. Reich, D. Shoikhet: Linearization models for parabolic 
dynamical systems via Abel's functional equation! Ann. Acad. Sci. Fen.! 35 2010 1-34

\art ELRS M. Elin, M. Levenshtein, S. Reich, D. Shoikhet: Some inequalities for the horosphere function and hyperbolically nonexpansive mappings on the Hilbert ball! J. Math. Sci. (N.Y.)! 201 2014 no. 5, 595-613 

\pre EJ M. Elin, F. Jacobzon: Parabolic type semigroups: asymptotics and order of contact!
Preprint, arXiv:1309.4002! 2013

\art ERS M. Elin, S. Reich, D. Shoikhet: A Julia-Carath\'eodory theorem for hyperbolically monotone mappings in the Hilbert ball! Israel J. Math.! 164 2008 397-411

\book ES  M. Elin, D. Shoikhet: Linearization models for complex dynamical systems. Topics in univalent functions, functional equations and semigroup theory! Operator Theory: Advances and Applications, 208. Linear Operators and Linear Systems Birkhäuser Verlag, Basel, 2010. xii+265 pp. ISBN: 978-3-0346-0508-3 

\art ESY M. Elin, D. Shoikhet, F. Yacobzon: Linearization models for parabolic type semigroups! J. Nonlinear Convex Anal.! 9 2008 205-214

\art F K. Fan: The angular derivative of an operator-valued analytic function! Pacific J. Math.! 121 1986 67-72

\art Fa P. Fatou: S\'eries trigonom\'etriques et s\'eries de Taylor! Acta Math.! 30 1906 335-400 

\art H M. Herv\'e: Quelques propri\'et\'es des applications analytiques d'une
boule \`a $m$ dimensions dans elle-m\^eme! J. Math. Pures Appl.! 42 1963 117-147

\book JP M. Jarnicki, P. Pflug: Invariant distances and metrics in complex analysis! Walter de Gruyter \& co., Berlin, 1993

\art Ju1 G. Julia: M\'emoire sur l'it\'eration des fonctions rationnelles! J. Math. Pures Appl.! 1 1918 47-245

\art Ju2 G. Julia: Extension nouvelle d'un lemme de Schwarz! Acta Math.! 42 1920 349-355

\art Ko A. Kor\'anyi: Harmonic functions on hermitian hyperbolic spaces! Trans.
Amer. Math. Soc.! 135 1969 507-516

\art K-S A. Kor\'anyi, E.M. Stein: Fatou's theorem for generalized
half-planes! Ann. Scuola Norm. Sup. Pisa! 22 1968 107-112

\art L-V E. Landau, G. Valiron: A deduction from Schwarz's lemma! J. London Math. Soc.! 4 1929 162-163

\art L L. Lempert: La m\'etrique de Kobayashi et la
repr\'esentation  des domaines sur la boule! Bull. Soc. Math. France! 109 1981 427-474

\art Li E. Lindel\"of: Sur un principe g\'en\'erale de l'analyse et ses
applications \`a la theorie de la repr\'esentation conforme! Acta Soc. Sci.
Fennicae! 46 1915 1-35

\art MM M. Mackey, P. Mellon: Angular derivatives on bounded symmetric domains! Israel J. Math.! 138 2003 291-315

\art N R. Nevanlinna: Remarques sur le lemme de Schwarz! C.R. Acad. Sci. Paris! 188 1929 1027-1029

\art RS1 S. Reich, D. Shoikhet: Semigroups and generators on convex domains with the hyperbolic metric! Atti Acc. Naz. Lincei Cl. Sc. Fis. Mat. Nat. Rend. Lincei! 8 1997 231-250

\book RS2 S. Reich, D. Shoikhet: Nonlinear semigroups, fixed points, and geometry of domains in Banach spaces! Imperial College Press, London, 2005

\book R W. Rudin: Function theory in the unit ball of $\C^n$! Springer, Berlin, 1980

\book S D. Shoikhet: Semigroups in geometrical function theory! Kluwer Academic Publishers, Dordrecht, 2001

\book St E.M. Stein: The boundary behavior of holomorphic functions of several
complex variables! Princeton University Press, Princeton, 1972

\art SW A. Sza􏰀\l owska, K. W\l odarczyk: Angular derivatives of holomorphic maps in
infinite dimensions! J. Math. Anal. Appl.! 204 1996 1-28

\art W S. Wachs: Sur quelques propri\'et\'es des transformations pseudo-conformes
avec un point fronti\`ere invariant! Bull. Soc. Math. Fr.! 68 1940 177-198

\art W\l1 K. W􏰀\l odarczyk: The Julia-Carath\'eodory theorem for distance decreasing
maps on infinite-dimensional hyperbolic spaces! Atti Accad. Naz. Lincei! 4
1993 171-179
 
\art W\l2 K. W􏰀\l odarczyk: Angular limits and derivatives for holomorphic maps of infinite dimensional bounded homogeneous domains! Atti Accad. Naz. Lincei! 5 1994 43-53

\art W\l3 K. W􏰀\l odarczyk: The existence of angular derivatives of holomorphic maps of
Siegel domains in a generalization of $C^*$-algebras! Atti Accad. Naz. Lincei! 5 1994 309-328

\art Wo J. Wolff: Sur une g\'en\'eralisation d'un th\'eor\`eme de Schwarz! C.R. Acad. Sci. Paris! 183 1926 500-502

\art Z J.M. Zhu: Angular derivatives of holomorphic maps in Hilbert spaces! J.
Math. (Wuhan)! 19 1999 304-308

\bye